\documentclass[11pt,twoside, reqno]{amsart}

\usepackage{amsmath}
\usepackage{amsthm}
\usepackage{amsfonts, amssymb}
\usepackage{mathrsfs}
\usepackage[all]{xy}
\usepackage{url}

\usepackage{latexsym}
\usepackage[dvips]{graphicx}

\theoremstyle{plain}
\newtheorem{lema}{Lemma}

\newtheorem{teo}[lema]{Theorem}

\theoremstyle{remark}

\newtheorem{obs}[lema]{Remark}

\theoremstyle{definition}
\newtheorem{defi}[lema]{Definition}

\newcommand{\Z}{\mathbb{Z}}
\newcommand{\R}{\mathbb{R}}

\renewcommand{\P}{\mathcal{P}}
\newcommand{\Q}{\mathcal{Q}}

\newcommand{\ZZ}{\Z[X,X^{-1},Y,Y^{-1}]}

\newcommand{\F}{\mathbb{F}_2}
\newcommand{\FF}{[\F,\F]}
\newcommand{\w}{\textrm{w}}

\begin{document}

\title[A counterexample to a strong version of the Andrews-Curtis conjecture]{A counterexample to a strong version of the Andrews-Curtis conjecture}

\author[J.A. Barmak]{Jonathan Ariel Barmak $^{\dagger}$}

\thanks{$^{\dagger}$ Researcher of CONICET. Partially supported by grant UBACyT 20020160100081BA}

\address{Universidad de Buenos Aires. Facultad de Ciencias Exactas y Naturales. Departamento de Matem\'atica. Buenos Aires, Argentina.}

\address{CONICET-Universidad de Buenos Aires. Instituto de Investigaciones Matem\'aticas Luis A. Santal\'o (IMAS). Buenos Aires, Argentina. }

\email{jbarmak@dm.uba.ar}

\begin{abstract}
\center\begin{minipage}{\dimexpr\paperwidth-8.4cm}
\medskip
\hspace{-0.7cm}{We prove that the presentations $\langle x,y | [x,y],1 \rangle$ and $\langle x,y | [x,[x,y^{-1}]]^2y[y^{-1},x]y^{-1},[x,[[y^{-1},x],x]] \rangle$}

\hspace{-0.7cm}{are not $Q^*$-equivalent even though their standard complexes have the same simple homotopy type.}
\end{minipage}
\end{abstract}



\maketitle

A finite presentation $\langle x_1,x_2,\ldots ,x_n | r_1,r_2,\ldots ,r_m\rangle$ of a group $G$ can be transformed into another presentation of $G$ by performing one of the following:

\noindent (i) Change a relator $r_j$ by $r_jr_i$ for some $i\neq j$.

\noindent (ii) Change a relator $r_j$ by $r_j^{-1}$.

\noindent (iii) Change a relator $r_j$ by a conjugate $gr_jg^{-1}$ for some $g$ in the free group $F(x_1,x_2,\ldots ,x_n)$.

\noindent (iv) Change each relator $r_j$ by $\phi (r_j)$ where $\phi$ is an automorphism of $F(x_1,x_2,\ldots ,x_n)$. 

\noindent (v) Add a generator $x_{n+1}$ and a relator $r_{n+1}$ which coincides with $x_{n+1}$.

\noindent (vi) The inverse of (v), when possible.

Moves (i) to (iii) are called $Q$-transformations, (i) to (iv) are $Q^*$-transformations and moves (i) to (vi) are called $Q^{**}$-transformations. Two finite presentations are said to be $Q$-equivalent ($Q^*$ or $Q^{**}$) if one can be obtained from the other by performing a sequence of $Q$-transformations ($Q^*$ or $Q^{**}$). The Andrews-Curtis conjecture states that any two balanced presentations ($m=n$) of the trivial group are $Q$-equivalent \cite{AC}. The weak version of the Andrews-Curtis conjecture states that two balanced presentations of the trivial group are always $Q^{**}$-equivalent. The latter is equivalent to the statement that any contractible finite $2$-dimensional CW-complex $3$-deforms to a point. The so called generalized Andrews-Curtis conjecture \cite[Section 4.1]{HM} says that any two presentations $\P$ and $\Q$ with simple homotopy equivalent standard complexes $K_{\P}$, $K_{\Q}$ are $Q^{**}$-equivalent. All these three conjectures are open. In this article we prove that there exist presentations $\P$ and $\Q$ with simple homotopy equivalent standard complexes which are not $Q^*$-equivalent. We do not know whether these presentations are $Q^{**}$-equivalent or not.

The main ingredients of the proof are an example of Evans \cite{Eva} of a matrix $M\in GL_2(\ZZ)$ which is not a product of elementary matrices and diagonal matrices, and a new invariant called the winding invariant. The present article was meant to be a section in the paper \cite{Bar} which introduces this invariant along with applications. We believe that it is better to present this example in a separate article.

Previously known methods for proving that presentations $\langle x_1,x_2,\ldots ,x_n | r_1,r_2,\ldots ,r_m\rangle$ and $\langle x_1,x_2,\ldots ,x_n | s_1,s_2,\ldots ,s_m\rangle$ are not $Q$ or $Q^{**}$-equivalent were developed by Browning \cite{Bro}, Hog-Angeloni and Metzler \cite[Section 2.2]{HM2} and Borovik, Lubotzky and Myasnikov \cite{BLM}. As explained in \cite[Section 2.2]{HM3}, the idea is to define a homomorphism $q$ from $F(x_1,x_2, \ldots, x_n)$ to a test group $G^*$ and prove that the $m$-tuples $(q(r_1),q(r_2),\ldots , q(r_m))$, $(q(s_1),q(s_2),\ldots , q(s_m))$ are not equivalent. The results of Hog-Angeloni and Metzler \cite[Theorems 2.3 and 2.4]{HM} show that solvable groups are not useful as test groups to distinguish $Q^{**}$-equivalence. In \cite{BLM, Bro2, HM3} there are results concerning the use of finite groups to distinguish $Q$ and $Q^{**}$-equivalences. Our methods appeared as a natural application when studying the winding invariant. However, they can be seen through this perspective as an application of solvable groups to distinguish $Q$-equivalence. Concretely we use $G^*=\F/\F''$, the free metabelian group of rank 2.  

The presentations $\P$, $\Q$ that we will define have fundamental group $\Z \times \Z$, so their Whitehead group is trivial. Therefore the torsion $\tau(f)\in Wh(K_{\P})$ of any homotopy equivalence $f:K_{\Q}\to K_{\P}$ is trivial. In order to distinguish $Q^*$-equivalence classes we will work in $GL_2(\Z[\pi_1(\P)])/GE_2(\Z[\pi_1(\P)])$ instead of $Wh(K_{\P})=GL(\Z[\pi_1(\P)])/GE(\Z[\pi_1(\P)])$. We recall the definitions of these concepts.

Given a ring $R$, denote by $E_n(R)$ the subgroup of $GL_n(R)$ generated by the elementary matrices. Recall that $E\in GL_n(R)$ is elementary if all the diagonal coefficients are $1\in R$ and all the other coefficients but one are $0\in R$. Note that $E_n(R)$ is a normal subgroup of $GL_n(R)$. We call $D_n(R)$ the subgroup of $GL_n(R)$ of diagonal matrices and $GE_n(R)$ the subgroup of $GL_n(R)$ generated by $E_n(R)$ and $D_n(R)$. If $R$ is a Euclidean ring, then $GE_n(R)=GL_n(R)$ for every $n\ge 1$. A ring $R$ is said to be generalized Euclidean if $GE_n(R)=GL_n(R)$ for every $n$. It was proved by Bachmuth and Moshizuki \cite[Theorem 1]{BM} that $R=\ZZ$ is not generalized Euclidean. Evans gives in \cite[Theorem C]{Eva} a concrete example of a $2$ by $2$ invertible matrix over that ring which is not in $E_2(R)$.

From now on $R$ will denote the ring $\ZZ$.

\begin{teo}[Evans] \label{evans}
The matrix 
\begin{displaymath}
\left(\begin{array}{cc}
1-2(X-1)Y^{-1} & 4Y^{-1} \\
-(X-1)^2Y^{-1} & 1+2(X-1)Y^{-1}
\end{array}\right)
\end{displaymath}
is in $GL_2(R)$ but not in $GE_2(R)$.

\end{teo}

We denote by $\F$ the free group on generators $x,y$. Its commutator subgroup is denoted $\FF$ or $\F'$. Recall that $w\in \F$ lies in $\FF$ if and only if the total exponents of $x$ and of $y$ in $w$ are both $0$. 

\begin{defi} \label{definicion}
Let $w\in \F$. Then $w=x_1^{\epsilon_1}x_2^{\epsilon_2}\ldots x_t^{\epsilon_t}$ where $x_i\in \{x,y\}$ and $\epsilon_i \in \{1,-1\}$ for each $i$. The word $w$ determines a path $\gamma_{w}$ in $\R^2$ which begins in $(0,0)$ and is a concatenation of paths $\gamma_1, \gamma_2,\ldots ,\gamma_l$. The path $\gamma_i$ moves one unit parallel to the axe $x_i$ and with positive or negative direction depending on the sign $\epsilon_i$. The image of $\gamma_w$ is contained in the grid $\Z\times \R \cup \R \times \Z$. Note that the ending point of $\gamma_w$ is $(k,l)$ where $k$ is the total exponent of $x$ in $w$ and $l$ is the total exponent of $y$. Suppose $w\in \FF$, so $\gamma_w$ finishes in $(0,0)$. For each $(i,j)\in \Z \times \Z$, let $a_{i,j}$ be the winding number $\w(\gamma_w, i+\frac{1}{2},j+\frac{1}{2})$ of $\gamma_w$ around the point $p=(i+\frac{1}{2},j+\frac{1}{2})$. Define the \textit{winding invariant} $P_w\in \Z[X,X^{-1},Y,Y^{-1}]$ of $w$ to be the Laurent polynomial $P_w= \sum a_{i,j} X^iY^j$. 
\end{defi}

Our notation for commutators is $[u,v]=uvu^{-1}v^{-1}$. So, for instance the winding invariant of $[x,y]$ is $P_{[x,y]}=1\in \ZZ$ and $P_{[y^{-1},x]}=Y^{-1}$. The path $\gamma_w$ is just the lift of $w$ to the Cayley graph $\Gamma(\Z \times \Z, \{x,y\})=\Z\times \R \cup \R \times \Z$, and the winding invariant $\F '\to \ZZ$ can be seen as the projection $N\to N / N'$ of the normal closure $N$ of $[x,y]$ in $\F$ onto the relation module $N/N'$ of the presentation $\langle x,y | [x,y] \rangle$. The geometric nature of our definition is useful to understand the intuition behind the algebraic arguments. The coefficients $a_{i,j}$ of the winding invariant appear in an article of Conway and Lagarias \cite[Section 5]{CL}. However, no use is given there to the ring structure of $\ZZ$ which is key in our results. 

The proof of the next result follows immediately from the definition or from the comments above about relation modules.

\begin{lema} \label{basica}
Let $w,w' \in \FF, u\in \F$. Then the following hold:

(i) $P_{ww'}=P_w+P_{w'}$.

(ii) $P_{w^{-1}}=-P_w$.

(iii) $P_{uwu^{-1}}=X^kY^lP_w$, where $k$ and $l$ are the total exponents of $x$ and $y$ in $u$.

(iv) $P_{[u,w]}=(X^kY^l-1)P_w$.
\end{lema}

Call a presentation $\P=\langle x,y | r_1,r_2,\ldots, r_m \rangle$ \textit{cocommutative} if each relator $r_j$ lies in $\FF$. 

\begin{obs} \label{gran}
Every presentation $Q^*$-equivalent to a cocommutative presentation is also cocommutative. We can associate with a cocommutative presentation $\P$ the vector $\Lambda(\P)=(P_{r_1},P_{r_2},\ldots, P_{r_m})\in \ZZ ^m$. The effect on $\Lambda(\P)$ of performing a  $Q$-trans\-for\-mation on $\P$ is to change a polynomial $P_{r_j}$ by $P_{r_j}+P_{r_i}$ for certain $i\neq j$, or by $-P_{r_j}$, or by $X^kY^lP_{r_j}$ for certain $k,l\in \Z$. Therefore, if $\P$ is a cocommutative presentation and $\Q$ is $Q$-equivalent to $\P$, then $\Lambda(\Q)^t=E\Lambda(\P)^t$ for some $E\in GE_m(R)$. Here $\Lambda(\P)^t, \Lambda(\Q)^t \in R^{m\times 1}$ denote the column vectors.
\end{obs}

If $\P=\langle x,y | r_1,r_2,\ldots ,r_m \rangle$ and $\Q=\langle x,y | s_1,s_2, \ldots, s_m \rangle$ are cocommutative presentations such that the normal closure of $\{r_1,r_2,\ldots ,r_m\}$ coincides with the normal closure of $\{s_1,s_2,\ldots, s_m\}$, then each $s_j$ is a product of conjugates of $r_i$'s and inverses of $r_i$'s. By Lemma \ref{basica}, $P_{s_j}$ is an $R$-linear combination of the $P_{r_i}$'s. Therefore $\Lambda(\Q)^t=M \Lambda (\P)^t$ for certain $M\in R^{m\times m}$. Symmetrically, $\Lambda (\P)^t=M' \Lambda (\Q)^t$ for some $M'\in R^{m\times m}$.

\bigskip

\noindent Let $$\P=\langle x,y | [x,y], 1 \rangle$$ and let $$\Q= \langle x,y | [x,[x,y^{-1}]]^2y[y^{-1},x]y^{-1}, [x,[[y^{-1},x],x]] \rangle .$$

\bigskip

\noindent We state the main result of the article.

\begin{teo}\label{main}
The standard complexes $K_{\P}$ and $K_{\Q}$ are simple homotopy equivalent while $\P$ and $\Q$ are not $Q^*$-equivalent.
\end{teo}

Since the Whitehead group of $\pi_1(K_{P})=\Z\times \Z$ is trivial \cite{BHS}, to prove simple homotopy equivalence we only need to prove homotopy equivalence. This can be achieved by standard methods and we postpone this part to the end of the proof.

\begin{figure}[tb] 
\begin{center}
\includegraphics[scale=0.75]{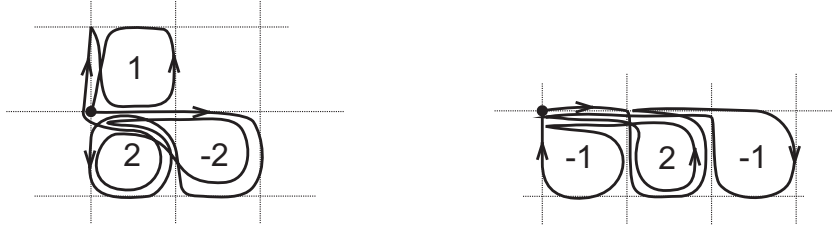}
\caption{At the left, the graphical representation of the curve $\gamma_{r_1}$ which begins in the black dot. The represented curve and the actual curve are homotopic in the plane with the centers of the squares removed. In the interior of the squares we see the corresponding winding numbers. At the right the curve $\gamma_{r_2}$.}\label{fig}
\end{center}
\end{figure}

We begin by proving that $\P$ and $\Q$ are not $Q$-equivalent. We compute first $\Lambda (\Q)$. Let $r_1=[x,[x,y^{-1}]]^2y[y^{-1},x]y^{-1}$ and $r_2=[x,[[y^{-1},x],x]]$ be the relators of $\Q$. By Lemma \ref{basica} the winding invariant of $r_1$ is $2(X-1)P_{[x,y^{-1}]}+YP_{[y^{-1},x]}$. Since $P_{[y^{-1},x]}=Y^{-1}$ and $P_{[x,y^{-1}]}=-P_{[y^{-1},x]}=-Y^{-1}$, then $P_{r_1}=1-2(X-1)Y^{-1}$ (See Figure \ref{fig}). Similarly $P_{r_2}=(X-1)P_{[[y^{-1},x],x]}=-(X-1)^2P_{[y^{-1},x]}=-(X-1)^2Y^{-1}$. Thus, $$\Lambda(\Q)=(1-2(X-1)Y^{-1}, -(X-1)^2Y^{-1})$$ is the first column of matrix $M$ in Theorem \ref{evans}.

On the other hand it is easy to see that $\Lambda (\P)=(1,0)\in R^2$. If $\P$ and $\Q$ are $Q$-equivalent, by Remark \ref{gran} there exists a matrix $E\in GE_2(R)$ such that $M \binom{1}{0}= \Lambda(\Q)^t=E\Lambda(\P)^t=E \binom{1}{0}$. Then $E^{-1}M \binom{1}{0}=\binom{1}{0}$. Thus $E^{-1}M$ is a matrix of the form 

\begin{displaymath}
\left(\begin{array}{cc}
1 & A \\
0 & B
\end{array}\right)
\end{displaymath} for some $A,B\in R$. Moreover, since $E, M\in SL_2(R)$, $B=1\in R$, so $E'E^{-1}M=Id$ for the elementary matrix 
\begin{displaymath}
E'=\left(\begin{array}{cc}
1 & -A \\
0 & 1 
\end{array}\right).
\end{displaymath}
Then $M=E(E')^{-1}\in GE_2(R)$, which contradicts Theorem \ref{evans}. This completes the proof that $\P$ and $\Q$ are not $Q$-equivalent. The last lines of our proof are implicit in the comments of \cite[p. 115]{DMV} about unimodular columns. The matrix $M$ was used by Myasnikov, Myasnikov and Shpilrain in \cite[Theorem 1.6]{Mya} to study $Q$-transformations of $m$-tuples in a non-free group $G$. As we mentioned above, the idea of our proof naturally appeared when studying applications of the winding invariant. We discovered later that our methods are very similar to those used in the proof of \cite[Proposition 5.1]{Mya}. 

The fact that $\P$ and $\Q$ are not $Q^*$-equivalent either will follow from the next lemma.

\begin{lema} \label{ene}
The normal closure $N$ of $\{r_1,r_2\}$ is $\FF$.
\end{lema}

Before we give a proof the lemma, we show how to use it to prove that $\P$ and $\Q$ are not $Q^*$-equivalent. Suppose they are. Then there is an automorphism $\phi$ of $\F$ such that $\phi \P=\langle x,y | \phi([x,y]), 1\rangle$ and $\Q$ are $Q$-equivalent. Since $Q$-transformations preserve the normal closure of the relators, the normal closure of $\phi([x,y])$ is $\FF$ by Lemma \ref{ene}. By a well-known result of Magnus, $\phi([x,y])$ is a conjugate of $[x,y]$ or $[x,y]^{-1}$. In any case, $\phi \P$ is $Q$-equivalent to $\P$, so $\P$ and $\Q$ are $Q$-equivalent, a contradiction.

Alternatively, that $\P$ and $\Q$ are not $Q^*$-equivalent follows from the fact that there is no matrix $E\in GE_2(R)$ such that $E\Lambda(\P)^t=\Lambda(Q)^t$ and the following. If $\phi \in \textrm{Aut}(\F)$, then the winding invariant $P_{\phi([x,y])}$ is a unit of $R$, so there exists $E'\in GE_2(R)$ such that $E' \Lambda(\P)^t = \Lambda(\phi \P)^t$ (see \cite{Bar}). Then $\phi \P$ and $\Q$ cannot be $Q$-equivalent.

\bigskip

\noindent\textit{Proof of Lemma \ref{ene}}. It is clear that $r_1,r_2\in \FF$, so we only need to show that $[x,y]=1$ in $\F/N$. Let $d=[x,y^{-1}]=xy^{-1}x^{-1}y$. Since $$1=r_2=[x,[d^{-1},x]]=xd^{-1}xdx^{-1}x^{-1}xd^{-1}x^{-1}d=xd^{-1}xdx^{-1}d^{-1}x^{-1}d$$ in $\F/N$, then $xdx^{-1}d^{-1}=dx^{-1}d^{-1}x=x^{-1}(xdx^{-1}d^{-1})x$. Therefore $e=[x,d]$ commutes with $x$ in $\F/N$.

\noindent On the other hand $1=r_1=[x,d]^2yd^{-1}y^{-1}=e^2yd^{-1}y^{-1}$, so 

\begin{equation}
d=y^{-1}e^2y.
\label{de}
\end{equation}

\noindent Finally, by definition $e=xdx^{-1}d^{-1}=x(y^{-1}e^2y)x^{-1}d^{-1}$. But since, $d=xy^{-1}x^{-1}y$, we have $dy^{-1}x=xy^{-1}$. We use this in the previous equation to obtain $$e=dy^{-1}xe^2x^{-1}yd^{-1}d^{-1}.$$

\noindent We use now that $e$ and $x$ commute to deduce $e=dy^{-1}e^2yd^{-2}=ddd^{-2}=1$ in $\F/N$. The last equality follows from (\ref{de}). By (\ref{de}) again we deduce, $d=1$ in $\F/N$. Thus, $d=[x,y^{-1}]\in N$, so $[x,y]\in N$.\hfill\ensuremath{\square}

\bigskip

To finish the proof of the theorem we need to prove that $K_{\P}$ and $K_{\Q}$ are homotopy equivalent. We have done the most important part already in Lemma \ref{ene}. It implies that $K_{\P}$ and $K_{\Q}$ have isomorphic fundamental groups with an isomorphism $f_*:\pi_1(K_{\Q})\to \pi_1(K_{\P})$ induced by a map which is the identity on $1$-skeletons. Then $K_{\P}$ and $K_{\Q}$ have isomorphic fundamental groups and the same Euler characteristic. In general this does not imply homotopy equivalence, but it does in our case since $\pi_1(K_{\P})$ is free abelian of rank $2$. This is explained by Harlander in \cite{Har}: Suppose $K$ and $L$ are finite connected two-dimensional complexes with $\pi_1(K)\simeq \pi_1(L)\simeq \Z\times \Z$ and $\chi(K)=\chi(L)$. Since $\Z\times \Z$ is aspherical, $H_2(\widetilde{K})$ and $H_2(\widetilde{L})$ are projective $R$-modules by the generalized Schanuel Lemma. By the Quillen-Suslin Theorem, these $R$-modules are free (see \cite{Swa} for Swan's comments on how to go from polynomials to Laurent polynomials). Since $\chi(K)=\chi(L)$, $H_2(\widetilde{K})$ and $H_2(\widetilde{L})$ have the same rank, so they are isomorphic. Once again, since $\Z \times \Z$ is aspherical, $H^3(\Z\times \Z, H_2(\widetilde{K}))=H^3(\Z\times \Z, H_2(\widetilde{L}))=0$, so $K$ and $L$ have isomorphic algebraic $2$-types. By MacLane-Whitehead's Theorem \cite[Theorem 4.9]{Sie}, $K$ and $L$ are homotopy equivalent. This finishes the proof of Theorem \ref{main}.


An alternative and simpler way to see that $K_{\P}$ and $K_{\Q}$ are homotopy equivalent is to construct a homotopy equivalence $K_{\Q}\to K_{\P}$ using standard techniques. Let $f:K_{\Q}\to K_{\P}$ be any map which is the identity on $1$-skeletons. This lifts to a map $\widetilde{f}:\widetilde{K}_{\Q}\to \widetilde{K}_{\P}$ between universal covers and we have the following commutative diagram of $R$-modules with exact rows

\begin{displaymath}
\xymatrix@C=30pt{ 0 \ar@{->}[r] & H_2(\widetilde{K}_{\Q}) \ar@{->}^{i}[r] \ar@{->}^{\widetilde{f}^{\sharp}_2}[d] & C_2(\widetilde{K}_{\Q}) \ar@{->}^{d_2}[r] \ar@{->}^{\widetilde{f}_2}[d] & C_1(\widetilde{K}_{\Q}) \ar@{->}^{d_1}[r] \ar@{=}^{\widetilde{f}_1}[d] & C_0(\widetilde{K}_{\Q}) \ar@{->}[r] \ar@{=}^{\widetilde{f}_0}[d] & \Z \ar@{->}[r] \ar@{=}[d] & 0 \\
                 0 \ar@{->}[r]  & H_2(\widetilde{K}_{\P}) \ar@{->}^{i'}[r]  & C_2(\widetilde{K}_{\P}) \ar@{->}^{d'_2}[r]  & C_1(\widetilde{K}_{\P}) \ar@{->}^{d'_1}[r]  & C_0(\widetilde{K}_{\P}) \ar@{->}[r]  & \Z \ar@{->}[r] & 0 }             
\end{displaymath}

\noindent A computation of Fox derivatives shows that $d_2$ has matrix representation

\begin{displaymath}
\left(\begin{array}{cc}
(1-y)(1-2(x-1)y^{-1}) & (x-1)^2(1-y^{-1}) \\
(x-1)(1-2(x-1)y^{-1}) & -(x-1)^3y^{-1} 
\end{array}\right)
\end{displaymath}
and $d'_2$ is represented by
\begin{displaymath}
\left(\begin{array}{cc}
1-y & 0 \\
x-1 & 0 
\end{array}\right).
\end{displaymath}

\noindent Note that $d_2=\binom{1-y}{x-1}\Lambda(\Q)$ and $d'_2=\binom{1-y}{x-1}\Lambda(\P)$. This is no coincidence as the boundary operator $C_2(\widetilde{K}_{\mathcal{T}})\to C_1(\widetilde{K}_{\mathcal{T}})$ is represented by $\binom{1-y}{x-1}\Lambda(\mathcal{T})$ for any cocommutative presentation $\mathcal{T}$ in which the normal closure of the relators is $\FF$.

Suppose $\widetilde{f}_2:C_2(\widetilde{K}_{\Q})\to C_2(\widetilde{K}_{\P})$ is represented by a matrix $\left(\begin{array}{cc}
a & b \\
c & d 
\end{array}\right).$
Then $d'_2\widetilde{f}_2=d_2$ implies that $a=1-2(x-1)y^{-1}$ and $b=-(x-1)^2y^{-1}$. The complex $K_{\P}$ is a wedge of a torus and a two-dimensional sphere, so $\widetilde{K}_{\P}$ is homeomorphic to the plane $\R^2$ with a sphere attached at every point with integer coordinates. It is easy to see that $H_2(\widetilde{K}_{\P})$ is a free $R$-module of rank $1$ and $i':H_2(\widetilde{K}_{\P})\to C_2(\widetilde{K}_{\P})$ is the inclusion in the second coordinate. Let $\gamma :C_2(\widetilde{K}_{\Q})\to H_2(\widetilde{K}_{\Q})$ be defined by the $1\times 2$ matrix $(4y^{-1}-c, 1+2(x-1)y^{-1}-d)$. There exist then a map $g:K_{\Q}\to K_{\P}$ which is the identity on $1$-skeletons and $\widetilde{g}_2=\widetilde{f}_2+i'\gamma :C_2(\widetilde{K}_{\Q})\to C_2(\widetilde{K}_{\P})$ (see \cite[Lemma 1.2]{Lat}). Therefore

\begin{displaymath}
\widetilde{g}_2=\left(\begin{array}{cc}
1-2(x-1)y^{-1} & -(x-1)^2y^{-1} \\
c & d 
\end{array}\right)+\left(\begin{array}{c}
0 \\
1 
\end{array}\right) \left(\begin{array}{cc}
4y^{-1}-c & 1+2(x-1)y^{-1}-d
\end{array}\right)=
\end{displaymath}
\begin{displaymath}
=\left(\begin{array}{cc}
1-2(x-1)y^{-1} & -(x-1)^2y^{-1} \\
4y^{-1} & 1+2(x-1)y^{-1}
\end{array}\right)
\end{displaymath} is the transpose of the matrix $M$ in Theorem \ref{evans}, which is invertible. Thus $\widetilde{g}_2:C_2(\widetilde{K}_{\Q})\to C_2(\widetilde{K}_{\P})$ is an isomorphism, and then so is $\widetilde{g}^{\sharp}_2 :H_2(\widetilde{K}_{\Q})\to H_2(\widetilde{K}_{\P})$. Therefore $g:K_{\Q}\to K_{\P}$ induces isomorphisms in $\pi_1$ and $\pi_2$, so it is a homotopy equivalence.

\bigskip

\noindent\textit{Comments.} The matrix $M^t\in GL_2(R)$ is the matrix used to compute the Whitehead torsion $\tau(g) \in Wh(K_{\P})$ of $g$ (see \cite[Section 2]{Lus}). Although $[M^t]=0\in Wh(K_{\P})$, we have used in the first part of the proof that $M\neq 0$ in $GL_2(R)/GE_2(R)$. 

Recall that by a result of Suslin \cite[Corollary 7.10]{Sus}, $GL_3(R)=GE_3(R)$. This seems to give some evidence that $\P$ and $\Q$ could be $Q^{**}$-equivalent. However, that $r,s\in \FF$ have the same winding invariant only says that $r$ and $s$ are equal modulo $\F''$ (see \cite{Bar}). The fact that $\P$ can be $Q^{**}$-deformed into a presentation $\P'$ with relators that are equal to those of $\Q$ modulo $\F''$ was already known by \cite[Theorem 2.3]{HM2}. Our results show that \cite[Theorem 2.3]{HM2} does not hold if we replace $Q^{**}$-transformations by $Q^*$-transformations: if $\P'=\langle x,y | s_1',s_2'\rangle$ and $\Q'=\langle x,y | r_1',r_2'\rangle$ are obtained from $\P$ and $\Q$ by $Q^*$-transformations, then $\Lambda(\P')^t=E_{\P}\Lambda(\P)^t$ and $\Lambda(\Q')^t=E_{\Q}\Lambda(\Q)^t$ for certain $E_{\P},E_{\Q}\in GE_2(R)$. On the other hand, if $r_js_j^{-1}\in \F''$ for $j=1,2$, then $\Lambda(\P')=\Lambda(\Q')$. This is not possible since $\Lambda(\P)$ and $\Lambda(\Q)$ are not in the same orbit under the action of $GE_2(R)$.


\begin{thebibliography}{99}

\bibitem{AC} J.J. Andrews, M.L. Curtis. \textit{Free groups and handlebodies}. Proc. Amer. Math. Soc. 16(1965), 192-195.

\bibitem{BM} S. Bachmuth, H.Y. Mochizuki. \textit{$E_2\neq SL_2$ for most Laurent polynomial rings}. Amer.
J. Math. 104(1982), 1181-1189.

\bibitem{Bar} J.A. Barmak. \textit{The winding invariant}. Preprint.

\bibitem{BHS} H. Bass, A. Heller, R. Swan. \textit{The Whitehead Group of a Polynomial Extension}.
 Publ. de L'Inst. des Hautes Etude Sci 22(1964).

\bibitem{BLM} A.V. Borovik, A. Lubotzky, A.G. Myasnikov. \textit{The finitary Andrews-Curtis conjecture, in: Infinite groups: geometric, combinatorial and dynamical aspects}. Progr. Math. 248(2005), 15-30.

\bibitem{Bro} W.J. Browning. \textit{The effect of Curtis-Andrews moves on Jacobian matrices of perfect groups}. Cornell University, Ithaca N.Y. (1976), unpublished.

\bibitem{Bro2} W.J. Browning. \textit{Normal generators of finite groups}. 1976, unpublished.

\bibitem{CL} J.H. Conway, J.C. Lagarias. \textit{Tiling with polyominoes and combinatorial group theory}. J. Combin. Theory Ser. A 53(1990), 183-208. 
 
\bibitem{DMV} K. Dennis, B. Magurn, L. Vaserstein. \textit{Generalized euclidean group rings}. J. Reine Angew. Math. 351(1984), 113-128. 
 
\bibitem{Eva} M. Evans. \textit{Primitive elements in the free metabelian group of rank 3}. J. Algebra 220(1999), 475-491.

\bibitem{Har} J. Harlander. \textit{A survey of recent progress on some problems in 2-dimensional topology, in: Advances in two-dimensional homotopy and combinatorial group theory}. London Math. Soc. Lect. Notes Series 446(2018), 1-26.

\bibitem{HM} C. Hog-Angeloni, W. Metzler. \textit{Geometric aspects of two dimensional complexes, in: Two-dimensional homotopy and combinatorial group theory}. London Math. Soc. Lect. Notes Series 197(1993), 1-35.

\bibitem{HM2} C. Hog-Angeloni, W. Metzler. \textit{The Andrews-Curtis conjecture and its generalizations, in: Two-dimensional homotopy and combinatorial group theory}. London Math. Soc. Lect. Notes Series 197(1993), 365-380.

\bibitem{HM3} C. Hog-Angeloni, W. Metzler. \textit{Further results concerning the Andrews-Curtis conjecture and its generalizations, in: Advances in two-dimensional homotopy and combinatorial group theory}. London Math. Soc. Lect. Notes Series 446(2018), 27-35.

\bibitem{Lat} M.P. Latiolais. \textit{Homotopy and homology classification of 2-complexes, in: Two-dimensional homotopy and combinatorial group theory}. London Math. Soc. Lect. Notes Series 197(1993), 97-124. 

\bibitem{Lus} M. Lustig. \textit{Nielsen Equivalence and Simple-Homotopy Type}. Proc. London Math. Soc. 62(1991), 537-562.

\bibitem{Mya} A.D. Myasnikov, A.G. Myasnikov, V. Shpilrain. \textit{On the Andrews-Curtis equivalence}. Contemp. Math., Amer. Math. Soc. 296(2002), 183-198.


\bibitem{Sie} A. Sieradski. \textit{Algebraic topology for two dimensional complexes, in: Two-dimensional homotopy and combinatorial group theory}. London Math. Soc. Lect. Notes Series 197(1993), 51-96.

\bibitem{Sus} A.A. Suslin. \textit{On the structure of the special linear groups over polynomial
rings}. Math. USSR Izv. 11(1977), 221-238.

\bibitem{Swa} R. G. Swan. \textit{Projective modules over Laurent polynomial rings}. Trans. Amer. Math. Soc. 237 (1978), 111-120.
        
\end{thebibliography}
\end{document}